\documentclass[10pt,aps,prl,twocolumn,nopacs,final,letterpaper,superscriptaddress,longbibliography]{revtex4-2}
\usepackage[utf8]{inputenc}
\usepackage{amsmath,amsfonts,amssymb}
\usepackage{graphicx}
\usepackage{array}
\usepackage{bm}
\usepackage{hyperref}
\usepackage{xcolor}
\usepackage{quotes}
\definecolor{goodblue}{RGB}{0, 91, 187}
\hypersetup{
  colorlinks=true,
  allcolors=goodblue,
  urlcolor=goodblue,
  citecolor=goodblue,
  pdfborder={0 0 0},
  breaklinks=true,
}

\begin{document}

\title{Exploring the limits of the law of mass action in the mean field description of epidemics on Erd\"os-R\'enyi networks}
\author{Francisco  J. Mu\~noz}
\email{francisco.munoz@urjc.es}
\affiliation{Departamento de Matem\'atica aplicada, ciencia e ingenier\'ia de materiales y tecnolog\'ia electr\'onica, Escuela Superior de Ciencias Experimentales y Tecnología, Universidad Rey Juan Carlos, 28933, M\'ostoles (Madrid), Spain}

\author{Luca Meacci}
\affiliation{ICMC, Universidade de Sao Paulo, Av. Trab. Sao Carlense 400, Sao Carlos, 13566-590, SP, Brazil}

\author{Juan Carlos Nu\~no}
\email{juancarlos.nuno@upm.es}
\affiliation{Departamento de Matem\' atica Aplicada, Universidad Polit\'ecnica de Madrid,  28040-Madrid, Spain}

\author{Mario Primicerio}
\affiliation{Dipartimento di Matematica "U. Dini", Universita degli Studi di Firenze, Viale Morgagni, Firenze, 28040, Italy}

\begin{abstract}
The manner epidemics occurs in a social network depends on various elements, with two of the most influential being the relationships among individuals in the population and the mechanism of transmission. In this paper, we assume that the social network has a homogeneous random topology of Erd\"os-R\'enyi type. Regarding the contagion process, we assume that the probability of infection is proportional to the proportion of infected neighbours.

We consider a constant population, whose individuals are the nodes of the social network, formed by two variable subpopulations: Susceptible and Infected (SI model). We simulate the epidemics on this random network and study whether the average dynamics can be described using a mean field approach in terms of Differential Equations, employing the law of mass action. We show that a macroscopic description could be applied for low average connectivity, adjusting the value of the contagion rate in a precise function. This dependence is illustrated by calculating the transient times for each connectivity. 

This study contributes valuable insights into the interplay between network connectivity, contagion dynamics, and the applicability of mean-field approximations. The delineation of critical thresholds and the distinctive behaviour at lower connectivity enable a deeper understanding of epidemic dynamics.
\end{abstract}

\maketitle

\section{Introduction}

Spreading of epidemics, or the transmission of any transmissible factor, either material or immaterial (e.g., innovation, fashion, political opinions, etc.), occurs through interactions between individuals. The society is represented as a fixed number $N$ of individuals linked in a social network, where individuals are nodes or vertices, and links or arcs signify contact or information exchange. We assume individuals are immobile, and the network structure remains constant over time. Individuals may belong to several networks based on the relationship considered, influencing disease transmission \cite{Granell}.

Mathematical modelling has proven insightful in understanding epidemic properties \cite{Kiss, Keeling1}. Two crucial factors in characterizing the spreading process are: (i) the structure of the population's network and (ii) the mechanism of interaction. Initial and boundary conditions also play a role in non-linear dynamics, where multiple equilibria could emerge \cite{Meacci} 

Regarding the network type, we assume spreading occurs in a classical Erd\"os-R\'enyi random network \cite{ER}. This undirected network consists of randomly linked nodes. For any node $i$, its connectivity $k_i$ is defined as the number of nodes linked to it, forming its neighbourhood. A limit case is a fully connected network (regular), where $k_i = k_{max} = N - 1$. In general, a probabilistic parameter $p \subset [0,1]$ is defined such that $p = k / k_{max}$, where $k$ is the average connectivity $k = <k_i>$ for $i=1,\ldots,N$. Here, $p$ represents the probability that two nodes are linked.

In the subsequent sections, we adopt epidemiological terminology, referring to "contagion" rather than "interaction" and "infective" rather than "influencer." We employ the simplest model, the so-called SI-model, where the population is divided into two subpopulations: infectives ($Z$) and susceptibles ($X$). The infection probability for a susceptible individual is linearly dependent on the number of infected individuals in its neighbourhood. This model assumes that an individual, once infected, remains so; there is no recovery or death from infection. This description falls within the class of agent-based models \cite{Bianchi2015}.

In epidemiology, a macroscopic description often uses the law of mass action \cite{Wilson, Voit}, postulating that the contagion's rate  is proportional to the product of the numbers of susceptibles $X(t)$ and infectives $Z(t)$. The constant of proportionality $\beta$ is assumed to be time-independent. The standard formulation assumes a homogeneously mixed population, where the functions $X(t)$ and $Z(t)$ coincide with the total average number of susceptible and infected individuals obtained from the agent-based model. However, this equivalence depends on the network structure. This description has been classically written in terms of Ordinary Differential Equations \cite{Bacaer}.

The paper's main goal is to demonstrate how the network's connectivity influences the mean field approach's adequacy in describing the time evolution of the two subpopulations. We show that a critical average connectivity $k_d$ exists, such that for $k > k_d$, the mean field approach provides an adequate description of the average number of infected individuals, and $\beta$ is independent of $k$. Conversely, for $k$ lower but still close to $k_d$, where a unique giant component exists, the mean field approach is useful only if $\beta$ is adjusted in terms of the average network connectivity. For lower connectivity, the existence of disconnected components prevents a global macroscopic description. Nevertheless, to study the effect of very low connectivity on the contagion parameter, we restrict the dynamics to the largest component and apply the mean field approach to obtain the dependence of $\beta$ on $k$. Additionally, these results are complemented with the study of the time the population takes to achieve equilibrium, allowing the detection of classical thresholds of Erd\"os-R\'enyi networks.

The remainder of the paper is organized as follows: Section 2 provides a detailed description of the network model and the macroscopic approach. Section 3 introduces the classical random Erd\"os-R\'enyi network (graph) and outlines its main structural properties in terms of the parameter $p$ or, equivalently, average connectivity $k$. This section also describes the main characteristics of the simulations performed to compare with the macroscopic description. In subsection 3.1, we calibrate the contagion parameter $\beta$ of the law of mass action in terms of $p$. The computation of transient and characteristic times in terms of $p$ is covered in the fourth section, where their dependence is depicted. Finally, we discuss these results and suggest how they can be translated to other more complex epidemic models.

\section{Dynamics on networks \textit{vs}. macroscopic models}\label{Sec2}

To start our simulations, we create the social network as an Erd\"os-R\'enyi model: $N$ individuals (nodes) are randomly generated and are connected by $L$ randomly placed links. If we define the connectivity of each node as the number of nodes that are connected to it, then the average connectivity is $2 \, L/N$. Equivalently, we can prescribe the probability $p$ that two nodes are linked which, as it can be easily shown, is given by: $p= 2 \, \frac{L}{N \, (N-1)}$.

With regards to the contagion, each individual can be susceptible ($X$) or Infected ($Z$). Since the size of the population is constant, then, $N = X(t) + Z(t)$ for all $t$.
Infectives do not change their state, while for any susceptible node the probability of contagion (more precisely the probability of becoming infected in a given time interval) depends on the fraction of infected that are in its neighbourhood, i.e. in the set of nodes connected to it. Formally,
 
\begin{equation}\label{Pcontagio}
P(susceptible \to infected) = \lambda \, \frac{\#infected}{\#neighbours} = \lambda \, \frac{Z_i(t)}{k_i}
\end{equation}
 where $Z_i(t)$ is the number of infected in the neighbourhood of node $i$ and $k_i$ is the connectivity of that node (that does not change with time). The parameter $\lambda\in [0, 1]$ is a measure of the efficiency of the contacts for the contagion to occur (and of the time interval chosen). 

In all the simulations, the size of the population, i.e. the number of nodes of the network, is fixed to $N = 10000$. This choice is the result of keeping a feasible computation time with a negligible internal noise. We run enough simulations, between 10 and 25, corresponding to the same initial number of infected nodes $Z_0$ randomly placed in the network. Averaging the number of infectives at each $t$, we obtain the function $Z(t)$. 

As in any other agent-based model, such as the cellular automata \cite{Wolfram, Gracia}, there is not a unique way of carrying out the simulations of the population dynamics on networks. Concerning the update of the state of the nodes, two main implementations are possible: (i) Synchronously, where all the nodes are updated at the same time step, using the information from the previous step, and recorded in the data file, and (ii) asynchronously, where, at each time step, each node is randomly selected and updated according to its current state and those of its neighbours and then, saved in the data file. It is important to note that the choice between these methodologies can impact simulation outcomes \cite{Porter}. 
Furthermore, there exists a straightforward relationship between the time scales in both methodologies: one time step in the synchronous simulation corresponds to $N$ time steps in the asynchronous case.

Although in this simple SI model, both synchronous and asynchronous implementations yield similar results, we opt for asynchronous simulations to better align with the mean field approach in the next section. The asynchronous simulations generate $N$ more points in a typical time evolution towards equilibrium, enhancing the accuracy of fitting the contagion parameter.

In the macroscopic model, our focus is not on the distribution of the two subpopulations in the network. Instead, we postulate a law of mass action, where the rate of increase of infectives (i.e., the rate of contagion) is proportional to the product of the total numbers, $X(t)$and $Z(t)$. This implies that the time evolution of the total numbers of the two subpopulations in the social network is described by the following Ordinary Differential Equations (ODEs):

\begin{eqnarray}\label{sysODEs}
X'(t) & = & - \frac{\beta}{N} \, X(t) \, Z(t)  \nonumber \\
Z'(t) & = & \frac{\beta}{N} X(t) \, Z(t) 
\end{eqnarray}
The parameter $\beta$ in the law of mass action describes the contagion process \cite{Wilson, Voit} and it is related but, not necessarily equal, to the parameter $\lambda$ of the network model. In any case, $\beta$ fixes the time scale of the dynamics and, in fact, enables a rescaling of the time to the dimensionless $\tau = \beta t$. 

Due to the total population constraint, this system can be reduced to only one ODE, the logistic equation:
\begin{eqnarray}\label{eqZ}
Z'(t) & = & \frac{\beta}{N} (N - Z(t)) \, Z(t) 
\end{eqnarray}
whose solution, for a given initial condition $Z(0) = Z_0$ is given by:
\begin{eqnarray}\label{soleqZ}
Z(t) & = & \frac{N \, Z_0} {Z_0+(N-Z_0)\, e^{-\beta \, t}} 
\end{eqnarray}
Note that the parameter $p$ of the Erd\"os-R\'enyi network does not appear explicitly in this system. As we have already mentioned, it is expected that for high connectivity (i.e. large $p$-values) and large network size $N$, this function approximates the dynamics of the number of infectives as obtained with the network model.


\section{The dynamics on a random Erd\"os-R\'enyi social network}\label{ERmodel}

\begin{figure*}
\centering
\includegraphics[width=0.95\textwidth]{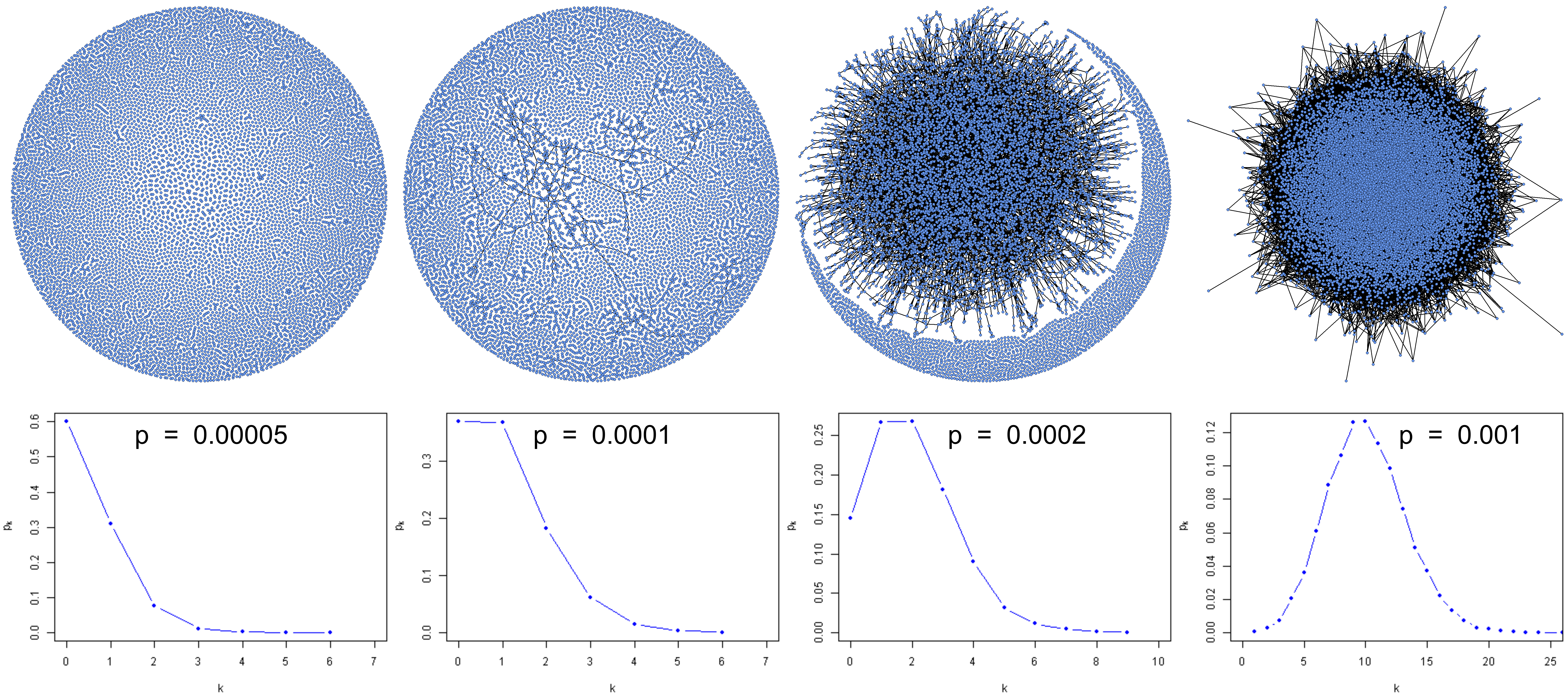}
\caption{Schematic representation of Erdős-Rényi networks is presented in the top row. The corresponding connectivity distributions are illustrated for four different values of $p$ in the bottom row. Notably, there is a discernible transition from a disconnected network to a fully connected network, with a Giant Component forming at $p=0.001$, just above the critical threshold $p_2$. This transition is also evident in the connectivity distributions: for $p$-values below the critical value $p_1$, the distributions lack a distinct peak, while for $p>p_2$, the connectivity distribution prominently peaks around its average connectivity $k$. In the case of $p=0.001$, this peak occurs at approximately $k = 10$. }
\label{kpmetrics}
\end{figure*}

No doubt, the structure of the network should be relevant for the dynamics of the contagion. In order to compare this dynamics with these network properties, let us delve into the details of the classical Erd\"os-R\'enyi network \cite{ER}. For each value of the parameter $p$, we can find the distribution of the connectivity $k_i$ of the nodes. It becomes apparent that for low $p$ values this distribution approximates a Poissonian (exponential) function, i.e. there is an average connectivity which coincides with the most frequent degree of the graph (see Fig.\ref{kpmetrics}).

\begin{figure*}
\centering
\includegraphics[width=0.9\textwidth]{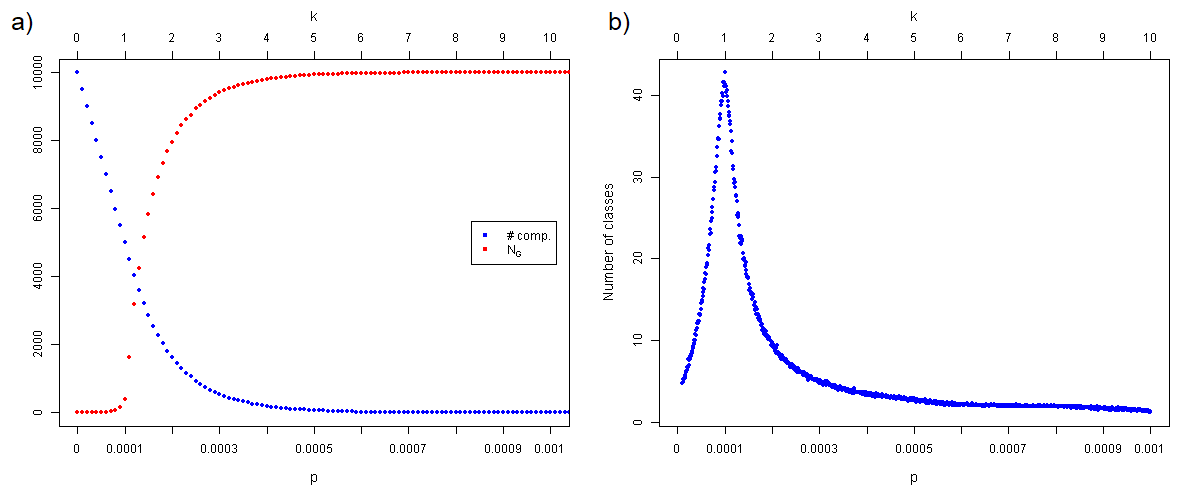}
\caption{(a) The number of components (blue) and the size of the largest component (red) are depicted as functions of $p$ (bottom $X$-axis) and the average connectivity $k$ (top $X$-axis) for an Erdős-Rényi graph with $N=10000$ nodes. Notably, the critical values $p_1 = 0.0001$ and $p_2 \approx 0.00092$, as discussed in the main text, mark the initiation of growth for the red curve and its approach to the graph size $N$, respectively. (b) The number of different sizes of components, averaged over 25 realizations, is presented as a function of $p$. Additional details, including explicit values for different component sizes, are provided in Table \ref{tb_p_classes_ng}.}
\label{NGdep}
\end{figure*}

For large $p$-values the network is formed by only one component (the so-called giant component), where all nodes can be reached from any other, meaning there is a path that connects any pair of nodes. In contrast, for low $p$-values the network disconnects, and several components appear. The transition between these phases occurs critically; two important $p$-thresholds are relevant to interpret the results we present in the next section (see Fig.\ref{NGdep}): (i) $p_1 = \frac{1}{N-1}$, estimating the $p$-value where the size of the largest component starts to increase exponentially as a function of $p$, and (ii) $p_2 = \frac{Ln(N)}{N}$, estimating the $p$-value where the size of the giant component attains the value of the network size, i.e. $N$. For $p > p_2$ there is a unique giant component containing most of the nodes of the network. Conversely, for $p < p_1$, the network is formed by low size disconnected components (see Fig. \ref{NGdep}b) to see the number of different sizes). This structural configuration prevents contagion from spreading between different components, affecting the entire network uniformly. With the network size set to: $N = 10000$, these thresholds take the values: $p_1 = 0.0001$ and $p_2\approx 0.00092$, which correspond to connectivity thresholds of $k_1 \approx 1$ and $k_2 \approx 9.2$, respectively.

As mentioned in the previous section, all simulations involve the generation of a random Erd\"os-R\'enyi network with parameter $p$ of size $N = 10000$. The probability of contagion of a susceptible node $i$ is calculated as in equation \ref{Pcontagio}. The initial populations are: $X_0 = N - Z_0$ and $Z_0$, whose individuals are randomly located in the nodes of the network. In order to diminish stochastic variability, each parameter setup is subjected to more than 10 simulations.

The dynamics of the SI model is straightforward for connected networks ($p > p_2$) where there is a path from each node to any other. In such cases, any initial number of infected individuals in the population leads to an epidemic, eventually infecting all susceptible individuals. Of course, complete contagion is unlikely if the network is disconnected unless each component initially contains at least one infected node.

To describe how the equilibrium situation of the system depends on $p$, we observe that:  

(i) For $p < p_1$, many isolated connected components of the network exist (some of them can simply be isolated nodes) (see Fig\ref{NGdep}). In this case, the components in which infective nodes are initially present will become totally infected, while the ones in which no infective nodes are initially present will be free of contagion forever. This implies that the dynamics of the system and its asymptotic equilibrium state, including the total number of nodes affected by the epidemic, depend not only on the initial number of infectives but also on their specific distribution within the network. 

This confirms the fact that for this situation a mean field approach in which the evolution of $Z(t)$ just depends on the total initial number of the infective individuals $Z(0) = Z_0$ cannot be found. 

(ii) For any value of p such that $p_1 < p < p_2$ there exit $m$ disconnected components an the number $m$ (average over several simulations) depends on $p$ (see Fig\ref{NGdep}a)). The same is true for the number $S$ of different sizes (Fig\ref{NGdep}b)). As $p$ increases, close to $p_2$, there is a largest component containing $N_G$ connected nodes and $m-1$ small components, disconnected from it and from each other, that contain $n_1, n_2, \ldots , n_{m-1}$ nodes. Of course, $N_G + \sum_1^{m-1} \, n_i = N$. If the initial situation involves infective individuals in the giant component, then the asymptotic number of the infected nodes will be $N_G + \sum n_j$, where the sum is extended to the disconnected components (if any) that initially contain infective nodes.

In the following table \ref{tb_p_classes_ng} and its corresponding figure \ref{NGdep}, we display the number of different sizes of the disconnected components averaged over 25 simulations for different values of $p$. For comparison, we also show in the table the values of $N_G$ and $m$ for these p-values.

\begin{table*}[!ht]
    \centering
    \begin{tabular}{|l|l|l|l|}
    \hline
        p & $N_G$ & $m$ & $S$ \\ \hline
        0.000010 & 4.76 & 9506.67 & 4.6 \\ \hline
        0.000050 & 17.88 & 7486.04 & 14.68 \\ \hline
        0.000075 & 49.16 & 6250.08 & 27.56 \\ \hline
        0.000100 & 403.68 & 5014.76 & 41.4 \\ \hline
        0.000125 & 3613.28 & 3842.08 & 27.44 \\ \hline
        0.000150 & 5796.8 & 2887.52 & 16.76 \\ \hline
        0.000200 & 7958.08 & 1625.88 & 9.4 \\ \hline
        0.000250 & 8941.48 & 918.88 & 6.32 \\ \hline
        0.000300 & 9406.24 & 542.6 & 4.8 \\ \hline
        0.000350 & 9657.96 & 322.88 & 4.08 \\ \hline
        0.000400 & 9803.72 & 188.92 & 3.36 \\ \hline
        0.000900 & 9998.68 & 2.32 & 1.68 \\ \hline
        0.001000 & 9999.72 & 1.28 & 1.28 \\ \hline
    \end{tabular}
\caption{The table presents the average number of different sizes $S$ of disconnected components, averaged over 25 simulations, for various values of $p$. Additionally, the table provides corresponding values for the size of the giant component $N_G$ and the average connectivity $m$ for each $p$-value. The network size is fixed at $N=10000$}
\label{tb_p_classes_ng}
\end{table*}

Denote by $N^*$ the total number of the nodes that belong to disconnected components that are initially infected. This means that the asymptotic number of infected individuals, $N_{eq}$ is given by: $N_{eq} = N_G + N^*$  (excluding the case in which the initial situation is such that no infected nodes are present in the giant component.)  Indeed, at equilibrium all the nodes of these components will be infected. Of course, $N^*$ increases with $Z_0$, because a larger number of initially infected nodes (i.e. $Z_0$) implies a higher probability that a disconnected component contains an infected node. Moreover, for a given $Z_0$, $N^*$ is decreasing w.r.t. $p$, because larger values of $p$ imply a smaller number of disconnected components. 

\begin{figure*}
\centering
\includegraphics[width=0.75\textwidth]{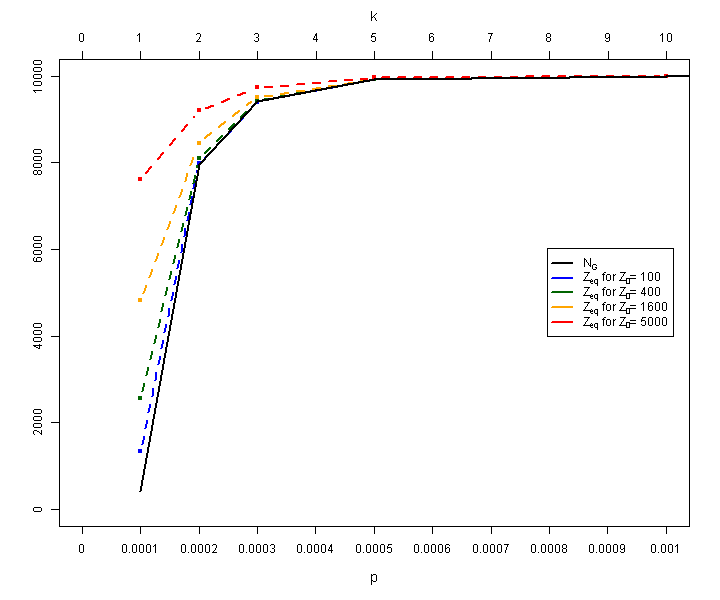}
\caption{Equilibrium values of the infected population obtained from simulations with varying initial populations, are depicted as a function of $p$. The results are compared with the size of the giant component. Note that for low values of $p$, $N_{eq}$ exceeds $ N_G$. However, as $p$ increases, the two values converge, with $N =N_{eq} = N_G$ holding true for $p > p_2$.}
\label{ZeqProbrule}
\end{figure*}

In Fig. \ref{ZeqProbrule} we display $N_{eq}$ as a function of $p$ for different values of $Z_0$. We see that in case $Z_0 = 100$, if $p > 0.0004$, then the corresponding (light blue) curve is close to the (black) curve representing $N_G$, i.e. $N^{*}$ ($= N_{eq} - N_G$) is negligible w.r.t. $N_G$. This means that for $p > 0.0004$ and "small" $Z_0$ epidemics predominantly takes place in the giant component and thus, we expect that the macroscopic description (with proper $\beta$) is possible. However, when $Z_0 = 2500$ (orange curve), the same situation occurs only for $p > 0.0005$ because this higher value of initial infective individuals implies that there are many more disconnected components in which the infection takes place. Hence $N^*$ is larger (for the same values of $p$). Thus, until $p$ is such that the number of disconnected components is small, a macroscopic description that holds from $t = 0$ is not feasible. 

(iii) For $p > p_2$ the giant component coincides with the whole network. Its connectivity increases with $p$, to reach complete (fully) connectivity for $p = 1$. In all cases, for $Z_0 > 0$, the whole population will become asymptotically infected. 
 
\section{Mean field calibration in terms of network connectivity}\label{calibration}

In general, the calibration of mean field models is not an easy problem \cite{Bicher}. Given an experimental model (microscopic), defined in terms of some variables and parameters that can be measured and controlled by the investigator and in terms of their (possibly stochastic) interactions, the question arises whether it is possible to describe its average behaviour using a mean field approach. In the affirmative case, a deterministic model can be written using differential equations (macroscopic) for the global variables and including parameters, that must be related to those used in the microscopic approach. The latter can also be thought of as a computational model that simulates the behaviour of the system from the bottom-up, i.e. by defining the rules that govern the dynamics. This is the situation that we are dealing with in section 2: the agent-based simulation of the contagion process that takes place in a social network. 

Two kinds of parameters can be controlled in the simulations: on the one hand, those related to the social network, specifically: (i) its size $N$ and (ii) the average connectivity $k$ or, equivalently, the probabilistic parameter $p$; on the other hand, the parameters concerning the population dynamics: (iii) the contagion rule and, once the probabilistic rule has been chosen, the contagion rate $\lambda$, (iv) the dimension of the initial populations $X_0 > 0$ and $Z_0 > 0$ ($X_0 + Z_0 = N$), and (v) the stopping time of the simulations, which it is automatic when the infected population achieves a stationary regime.

The mean field model \ref{eqZ} includes the parameter $\beta$ that has to be related to the corresponding parameter $\lambda$ of the agent-based model. The contagion rule applied to each susceptible individual depends on the proportion of infected in its neighbourhood, but also on the efficiency of the contacts measured by $\lambda$. The interesting question arises about how $\beta$ differs from $\lambda$ in terms of the local connectivity, once the size $N$ of the network is fixed. The strategy applied for the calibration consists in comparing the time evolution of the populations obtained by the agent-based model and by deterministic model \ref{eqZ} (of course, more sophisticated and rigorous procedures exist \cite{Bicher}). 
 
To do this comparison, it will be useful to introduce the following dimensionless variable $R(t)$, representing the ratio between the susceptible and the infected individuals at any time $t$: 
\begin{equation} 
R(t) =  \frac{X(t)}{Z(t)}                                      \end{equation}
It is immediately seen that $R(t)$ satisfies a linear ODE. The derivative of the previous expression yields:
\begin{equation} 
R'(t) = \frac{X'(t)}{Z(t)} - \frac{X(t)\,Z'(t)}{Z^2(t)} = \frac{\beta}{N} \frac{X(t)}{Z(t)} \, (X(t) + Z(t)) 
\end{equation}
and, since $X(t) + Z(t) = N$ for all $t$, then, we obtain the ODE:
\begin{equation} 
R'(t) = - \beta  \, R(t)                                       \end{equation}
The solution of the Initial Value Problem with $R(0) = R_0$ is given by: 
\begin{equation}\label{RODE}  
R(t) = R_0 \, e^{ - \beta \, t}                                \end{equation}

We next study the problem of $\beta$ calibration as a function of the linking probability $p$:

Case (i): For a population forming a fully connected network ($p = 1$) of size $N$, the parameters $\lambda$ and $\beta$ coincide. Indeed, for each susceptible individual $X$ the probability of changing its state per unit time, is given by $\lambda$ multiplied by the fraction $Z/N$ of the infected individuals in its neighbourhood that, in this case, coincides with the whole network. Hence, summing up over all $X$, the rate of change of the number of the infected nodes is: 
\begin{equation}\label{lawM}
P(t) = \frac{\lambda}{N} \, X(t) \, Z(t)
\end{equation}
which is the ODE corresponding to a macroscopic law of mass action with $\lambda = \beta$. 

We can easily show that the results of the simulations confirm this fact. For each simulation, we can plot the logarithm, $\ln{\left(\frac{X(t)}{Z(t)} \right)}$, where $X(t)$ and $Z(t)$ are the susceptible and infected population for each $t$. The graph of this logarithm is a straight line with slope $m = -\beta$ and intercept $n = \ln(R_0)$. Thus, we can apply a standard linear regression and find the value of $\beta$ for any simulation (with given initial condition $Z_0$). These $\beta$-values are averaged over 25 simulations to obtain the average $\beta$ for each $p$. As expected, the result is $\beta= 1$ with a standard error of the order of $10^{-4}$,  irrespectively of the value of $Z_0$.

Naturally, the expected error decreases with the number of simulations. Additionally, it becomes evident that it decreases when the network size $N$ increases. This reduction in error is attributed to the fact that the values obtained by the simulations are integers while the solutions of \ref{eqZ} are real numbers. Notably, this "rounding" error decreases with a larger network dimension (see table \ref{p1SR}).

\begin{table*}[!ht]
    \centering

    \begin{tabular}{|l|l|l|l|l|}
    \hline
        N & $Z_0$ & $\beta$ & Std. Error  \\ \hline
            & 10 &1.016& 6.695E-03 \\  
        100 & 25 &0.955&9.273E-03\\ 
            & 50 &1.069&1.669E-02\\ \hline
            
            & 10 &1.024& 5.438E-04 \\  
        1000 & 50 &0.981&8.242E-04\\ 
            & 500 &0.974&1.961E-03\\ \hline
            
            & 10 &0.997& 1.193E-04 \\  
        5000 & 250 &1.005&2.664E-04\\ 
            & 2500 &1.022&5.284E-04\\ \hline
            
            & 10 &0.995&6.073E-05 \\  
        10000 & 500 &1.004&1.251E-04\\ 
            & 5000 &0.998&2.306E-04\\ \hline
    \end{tabular}
    \caption{Results of the estimation of the contagion rate $\beta$ for the case $p=1$ are presented for various values of $N$ and $Z_0$. The estimations are based on the average of 25 simulations.}
    \label{p1SR}
\end{table*}

\begin{figure*}\
\centering
\includegraphics[width=0.9\textwidth]{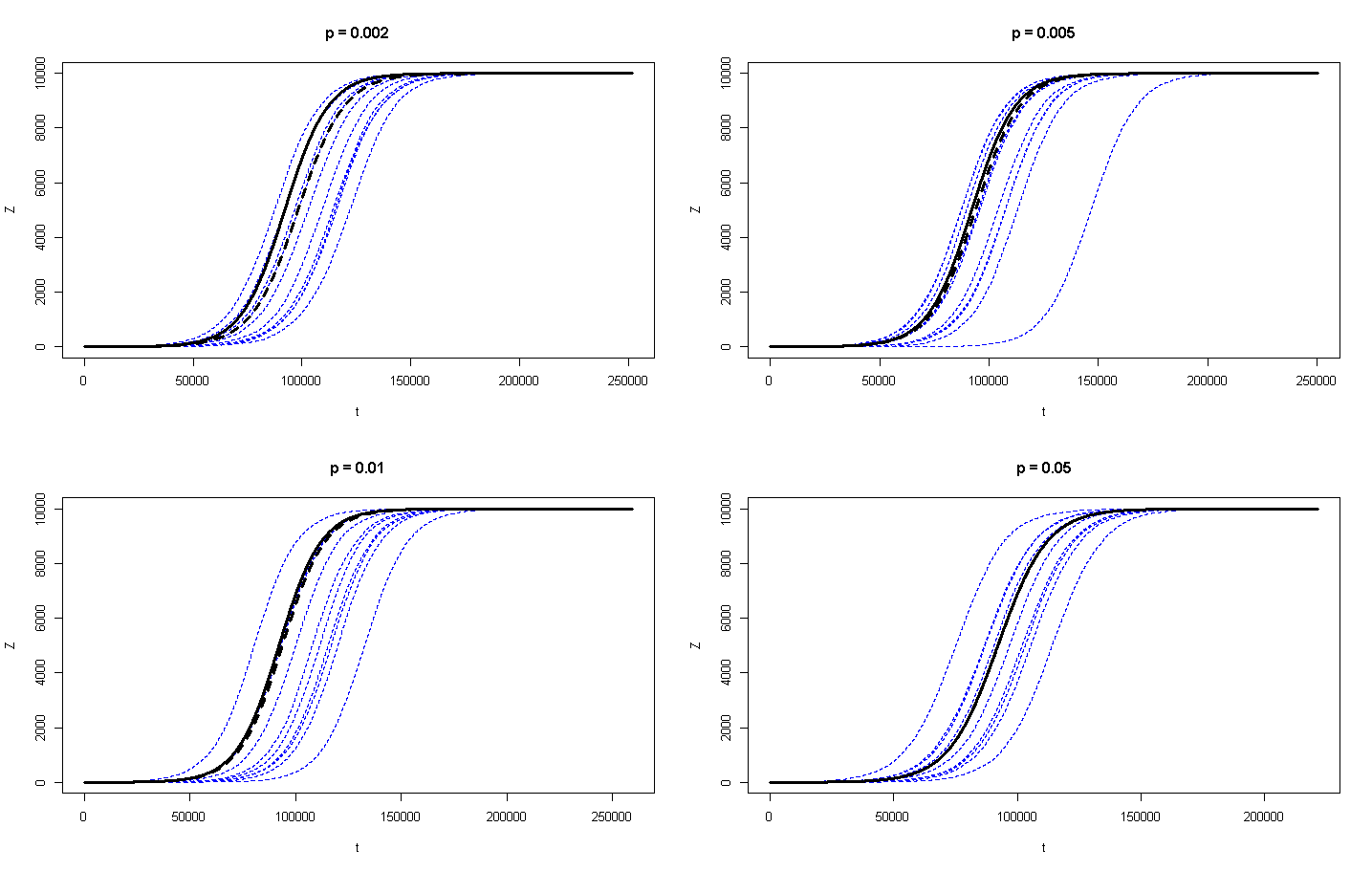}
\caption{The time evolution of the infected population $Z$, originating from a single infective individual randomly located in an Erdős-Rényi network, is depicted for four different values of $p$: $p=0.002, 0.005, 0.01, 0.06$. Each panel shows results from 10 simulations (blue), the average curve from these simulations (dashed black) and the analytical curve (solid black) obtained from equation \ref{eqZ} with $\beta=1$. Notably, with increasing $p$, the alignment between the blue and black curves becomes nearly complete.}
\label{z01}
\end{figure*}

Case (ii): When the network forms a unique giant component, i.e. $N_G=N$, it could be expected that the situation coincides with the case of fully connected network. However a simple experiment reveals that this is not the case, at least for all $p$ belonging to the interval $(p_2,1)$. To illustrate this, consider the limit case $Z_0=1$ and conduct simulations for different values of $p$. In Figure \ref{z01}, the curves representing $Z(t)$ for $p=0.002$, $p=0.005$, $p=0.01$ and $p=0.05$ are displayed. Each graph showcases the curves corresponding to 10 different simulations.


\begin{figure*}\
\centering
\includegraphics[width=0.95\textwidth]{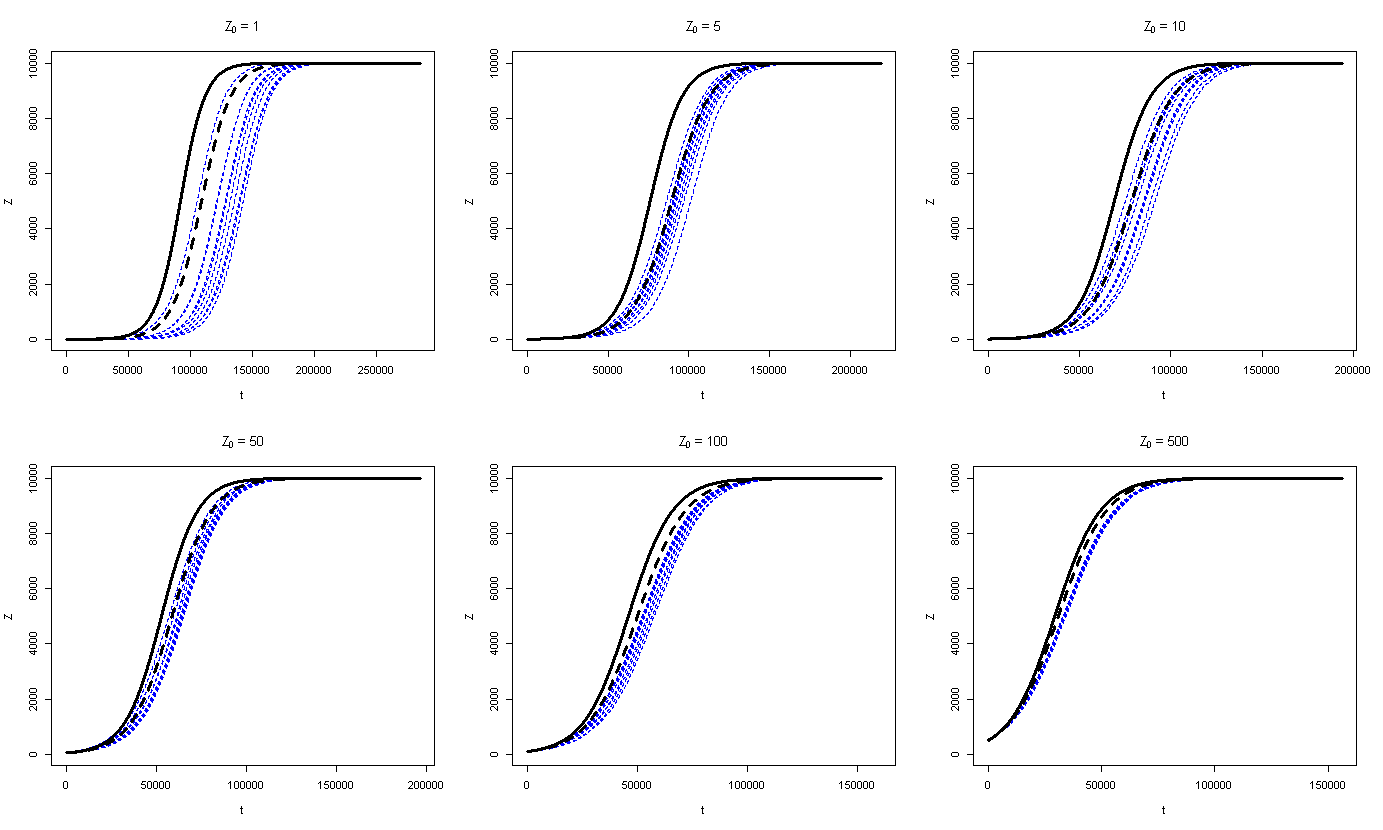}
\caption{The time evolution of the infected population for $p=0.001$ is presented for various initial infective populations: $Z_0 = 1, 5, 10, 50, 100, 500$. The analytical solution of the Initial Value Problem associated with equation \ref{eqZ} is represented by the solid black curve. 
Increasing $Z_0$ results in a decrease in the standard deviation $SD$ of experimental curves, as a larger sample of node degrees is considered. Nonetheless, despite this effect, the discrepancy between the average curve (dashed black) and the analytical curve persists due to the low average connectivity of the network for this $p$-value.}
\label{p0.001}
\end{figure*}

It is evident from figures \ref{z01}c and \ref{z01}d) that for larger $p$-values the 10 curves tend to coincide on average with the theoretical curve with $\beta=1$ (depicted in black). Conversely, for smaller values of $p$ (figures \ref{z01}a and \ref{z01}b) the curves are more distant. This discrepancy indicates that the influence of weakly connected nodes is more pronounced for smaller $p$ values, leading to a more significant dispersion of the $Z=Z(t)$ curves; when $Z_0=1$ the unique infected individual can reside in a node with a connectivity from the left "tail" of the degree distribution (see Fig. \ref{kpmetrics}). The dispersion in the curves of figures \ref{z01} can be quantified by the Standard Deviation $SD(p)$ of the values $\tau$ such that $Z(\tau)=N/2$.  As illustrated in Fig. \ref{p0.001} for $p=0.001$, $SD(p)$ decreases as $Z_0$ values increase.  As we will show below, this is a consequence of the lower average connectivity of the network for low $p$-values, which slightly slows down the dynamics.


The observed dependence of the dynamics on the initial infected population, even for $p$-values that yield a unique giant component, prevents the application of a mean field approach based on the ODE \ref{eqZ}. This prompts the question whether there exist a critical $p$-value above which the mean field approach, relying on the law of mass action, becomes statistically feasible. Given that the most unfavourable scenario occurs when $z(0)=1$, i.e. starting from a unique infected individual, we conducted 25 simulations for various $p$ values. In each simulation, the infected individual was placed randomly in a node of a random generated Erd\"os-R\'enyi network.


\begin{figure*}
\centering
\includegraphics[width=0.7\textwidth]{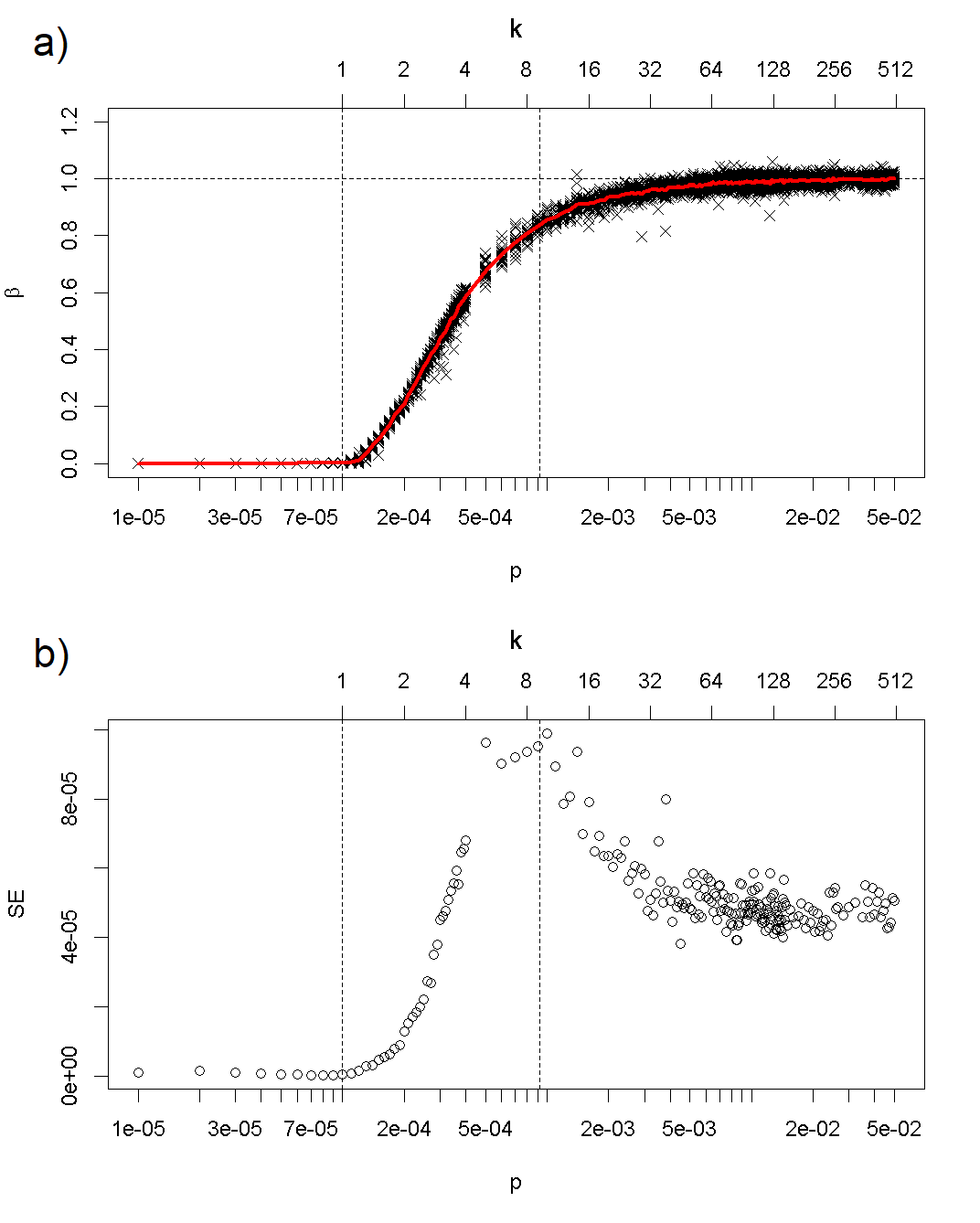}
\caption{(a) Semi-log plot illustrating the estimated $\beta$-values as a function of the linking probability $p$. As it can be observed, $\beta$ exhibits a slow increase from $p=0$ ($\beta(0)=0$) to $p=p_1$ ($\beta(p_1) \approx 10^{-3}$), followed by a sharp rise until it asymptotically converges to 1. As discussed in the main text, $\beta$ statistically equals 1 for $p >p_d=0.02$, approximately an order magnitude larger than the percolation threshold $p_2=0.00092$. The red curve represents the average $\beta$-values for each $p$.
(b)  Semi-log graph illustrating the averaged Standard Error (SE) in the estimation of $\beta$ as a function of $p$. Notably, the SE reaches its maximum value precisely at $p=p_1$, where the number of sizes in the network fragmentation also exhibits a peak (see Fig. \ref{NGdep}b). For larger $p$-values, the averaged SE tends asymptotically to $\approx 4 \, 10^{-5}$. }
\label{panel}
\end{figure*}


For each $p$-value, we perform a $t$-test over the distribution of estimated $\beta$ obtained in 25 simulations to assess whether the average is statistically equal to 1 within a confidence level of $95\%$. As required, the hypothesis of normality of the sample was verified. Exploring the range $(p_2,1)$, we identify a threshold, $p_d$ that separates the $p$-values into two sets: one formed by $p$-values that satisfies the test, i.e. $\beta$ is statistically equal to unity for them, and another formed by $p$-values for which $\beta$ differs from 1 statistically. For $N = 10^4$, this critical value is $p_d \approx 0.049$.  Consequently, $\beta =1$ with an Standard Error of the order of $10^{-5}$ for $p > p_d$ (see Fig. \ref{panel}(b)). It is worth noting that $p_d$ is an order of magnitude larger than the structural threshold $p_2 \approx 0.00092$ (see Fig. \ref{panel}(a)).

\begin{figure*}
\centering
\includegraphics[width=0.7\textwidth]{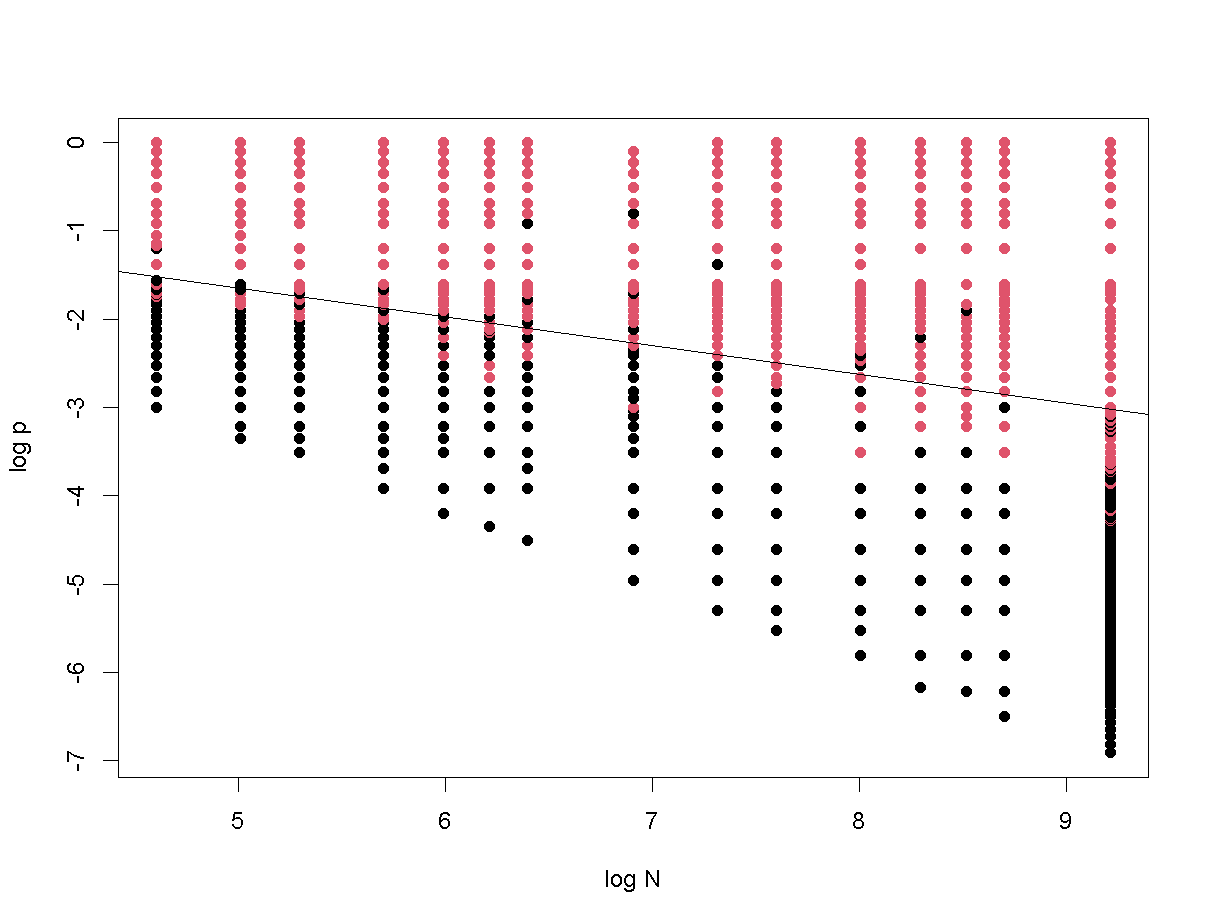}
\caption{Log-log representation of the set of $p$-values for each network size $N$ with indications of whether they satisfy the $t$-test for the corresponding $\beta$-value. Points in red denote cases where $\beta$ is statistically equal to 1, while black points indicate cases where it is not. The application of SVM algorithm yields a separatrix line in the $(\log(N),\log(p))$ plane, defining the threshold $p_d$. Remarkably, $p_d$ seems to follow a power law of $N$, specifically, $p_d \approx N^{-0.33}$. The lowest black points for each $N$ are close to $p_2$.}
\label{svm}
\end{figure*}

As expected, the threshold $p_d$ depends on the network size $N$. Figure \ref{svm} illustrates the results obtained from the $t$-test applied to the 25 simulations across the same range of $p$-values for different network sizes $N$. The red points correspond to the $p$-values for which $\beta=1$ (within a $95\%$ confidence level), whereas the black points represent $p$-values where $\beta < 1$. To define the threshold $p_d$ we apply a Support Vector Machine (SVM) algorithm to this set of points \cite{Meyer}. As shown, there is a separatrix line that separates both subsets (red-black) and it can be described by the log-log equation:
\begin{equation}
\log(p_d) = 0.03 - 0.33  \, \log(N)
\end{equation}
that yields the power law:
\begin{equation}
p_d \sim N^{-0.33}
\end{equation}


Case (iii): We now turn our attention to situations where the network is fragmented into more than one component. As depicted in Figure \ref{NGdep}, the number of components decreases with $p$ until reaching the percolation threshold at $p=p_2$. 

There are two factors that prevents a mean field approach in this $p$-interval: (i) the dependence with the initial conditions caused by the existence of several components and (ii) the low averaged connectivity that decreases the local contagion rate and, consequently, the value of $\beta$. 

To address the dependence on the initial conditions, we constrain the dynamics to the largest component for each $p$. This enables us to study the dependence of $\beta$ on $p$, under the assumption that a mean field approach is feasible for this component. Let $x(t)$ and $z(t)$ represent the numbers of the susceptible and infectious individuals in that component for a given $p$-value. Since the giant component has no connections with the rest of the network, $x(t) + z(t) = N_G$ and we can tentatively write the mean-field ODE for the giant component in the form: 
\begin{equation}\label{eq5} 
z'(t)   =   \frac{\beta}{N_G} (N_G-z(t)) \, z(t) 
\end{equation}
Consequently, we can define $r(t)=x(t)/z(t)$ and apply linear regression to the points representing the values of $ln(r(t))$ as obtained by the simulations, focusing on counting the infectives in the giant components.

Applying the same procedure, using a unique individual as a probe to explore the SI dynamics on the network, the logistic curves of the simulations are employed to derive an averaged $\beta$-value and its corresponding Standard Deviation (see Fig. \ref{panel} B). Notably, from an almost null value for $p < p_1$, $\beta$ increases continuously with $p$ until asymptotically reaching its macroscopic value of 1. It's worth mentioning that the function $\beta=\beta(p)$ is precisely defined at the extreme values of the $p$-interval: $\beta(0) = 0$ and $\beta (1) = 1$. 

The increase of the value of $\beta$ in the interval $(p_1,p_2)$ is strictly due to the rise in the average connectivity of the network. By using this unique initial infected individual within the largest component (not necessarily giant, as shown in Fig.\ref{panel}), we are neglecting the effect of the disconnected components in the population dynamics which, as we discussed earlier, is reflected in the dependence of the dynamics on the initial conditions. 

It is noteworthy to emphasize the behaviour of the averaged Standard Error (SE) of $\beta$ (panel (b) in figure \ref{panel}). In contrast to the averaged $\beta$-value that, as commented above, increases monotonously, the function SE(p) reaches a maximum at $p = p_1$. This critical $p$-value is defined based on the structure of the E-R networks as the inverse of network size. As observed before, it is attained at the same $p$-value where the distribution of sizes of disconnected components is maximum (see Fig. \ref{panel}). Here, we detect this threshold as a property of the dynamics of a population that follows a SI model with a probabilistic local contagion.

\section{Transient and characteristic times}

The results derived from the preceding section should yield a straightforward impact on the duration it takes for epidemics to affect the entire population. This duration can be gauged through simulations, specifically by assessing the average number of time steps required to reach population equilibrium, signifying complete infection. Alternatively, the analytical solution of the Ordinary Differential Equation described in Equation \ref{eqZ} can be employed to measure the time evolution of the average infected population. 

As demonstrated earlier, this latter approach holds validity only for $p >p_d$. In such instances, the time scale of equation \ref{soleqZ} is appropriately measured by $\tau = \frac{1}{\beta}$, regardless of the values of $p$, the network size $N$ or the initial population $z_0$. Unlike the reciprocal of the largest eigenvalue, the characteristic time takes into account the entire trajectory leading to the limiting value and, in so doing, it is able to differentiate systems with the same largest eigenvalue but exhibiting distinct dynamics \cite{Llorens}. This characteristic time serves as a metric for the duration a system requires to attain asymptotic behaviour, such as equilibrium. Naturally, the characteristic time is contingent upon the values of $N$ and $z_0$. 

For the ODE \ref{eqZ}, the characteristic time can be defined as follows:
\begin{equation}
T_c = \frac{\int_0^{\infty} t \, Z'(t) \, dt}{\int_0^{\infty} Z'(t) \, dt}
\end{equation}
where $Z(t)$ is the explicit solution of \ref{eqZ}, given by formula \ref{soleqZ}. A simple calculation yields:
\begin{equation}
T_c = \frac{N}{N - Z_0} \, \ln{\left(\frac{N}{Z_0}\right)} \, \frac{1}{\beta}
\end{equation}
independently of $p$. Specifically, for $\beta = 1$ and $Z_0 = 1$ this expression reduces to:
\begin{equation}
T_c = \frac{N \, \ln{\left( N \right)}}{N - 1}
\end{equation}
that, for the particular size $N =10^4$, particularizes to $T_c \approx 9.21$.

It is important to note that this time can be calculated whenever a solution to an ODE is accessible. Therefore, it remains applicable even when dealing with scenarios where $p < p_d$, limiting the dynamics to the largest connected component. In such cases, it is essential to consider the component size $N_G$ instead of the total network size $N$, along with the associated values of $\beta$ corresponding to each $p$ (as discussed in the previous section).

Another measures of duration applicable to discrete systems are the so-called $t_{\rho}$, representing the time the population takes to achieve a specified percentage $\rho$ of the equilibrium population (as discussed in \cite{Llorens} and references therein).  These times can be directly inferred from simulations, pinpointing the initial instance when the population, starting at $Z_0$, surpasses the defined threshold of $L = \rho , N_{eq}$.

\begin{figure*}
\centering
\includegraphics[width=0.8\textwidth]{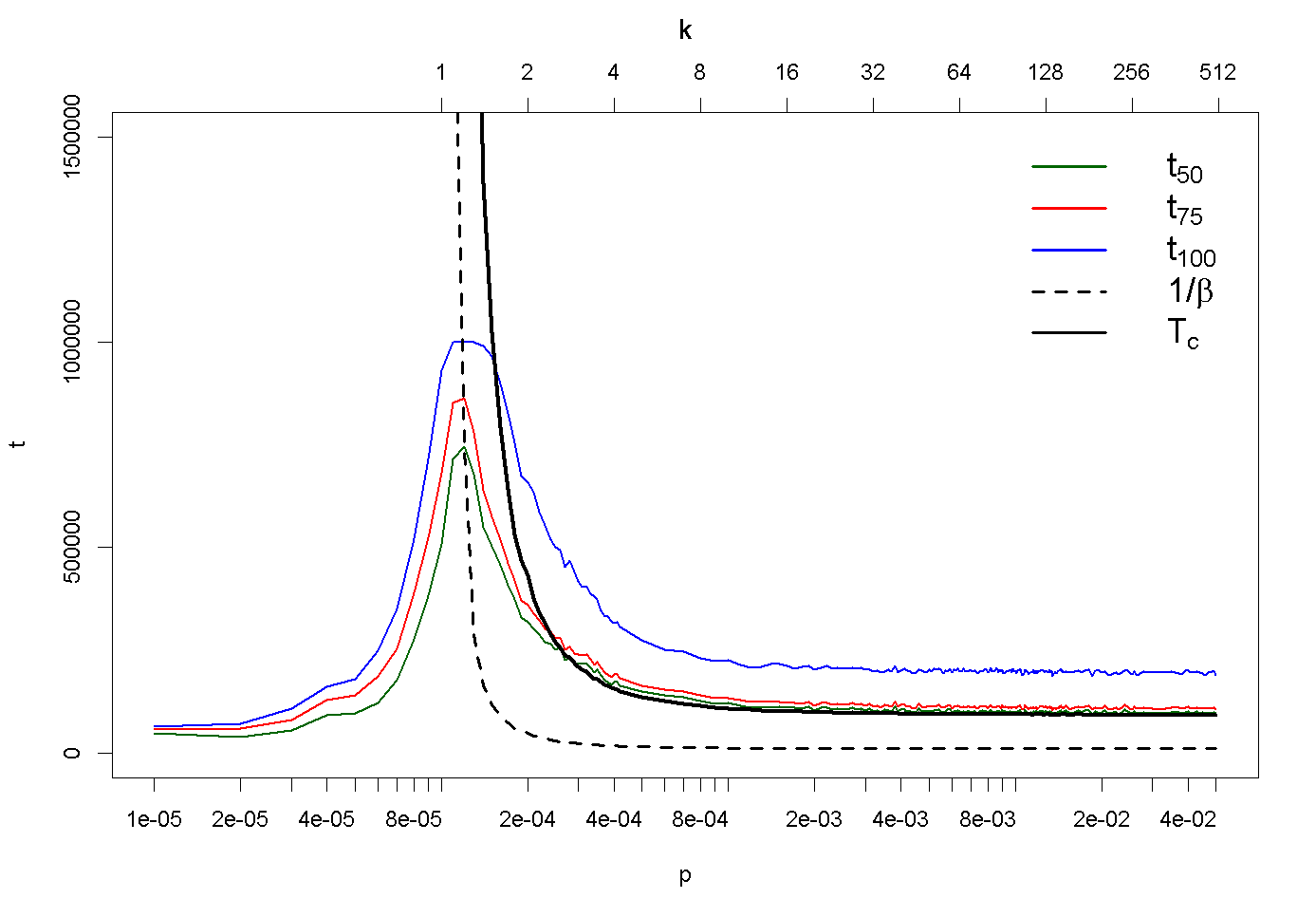}
\caption{Semi-log plot of the transient times as a function of $p$ (bottom axis) and $k$ (top axis). The times $t_{\rho}$ for $\rho = 50, 75$ and $100$ is represented along with the characteristic time $T_c$ and the reciprocal of the eigenvalue associated with the ODE \ref{eqZ}, i.e $1/\beta$. For $p > 10^{-3}$, all three times exhibit similar behaviour, remaining constant with $t_{\rho} > T_c > \frac{1}{\beta} = 1$. However, while the times $t_{\rho}$ reach a peak at $p \approx 10^{-4}$, both $T_c$ and $1/\beta$ monotonously decrease for all $p$ (approaching $\infty$ as $p \to 0$). The $Y$-axis represents simulation time in asynchronous updating.}
\label{tt}
\end{figure*}

Similar to the characteristic time, when dealing with scenarios where $p < p_2$, we calculate $t_{\rho}$ from simulations confined to the largest component, as illustrated in Fig. \ref{tt}. The transient times $t_{50}$, $t_{75}$ and $t_{100}$ are visually presented in Fig. \ref{tt}. It is evident that they follow the order $t_{50} < t_{75} < t_{100}$ and all four exceed both $T_c$ and $\beta^{-1}$. It is noteworthy that for $p > p_2$, these times remain relatively constant, indicating that substantial changes in network connectivity have minimal impact on the duration required for the complete spread of the epidemic.

As evident in Fig. \ref{tt}, the intriguing behaviour becomes apparent for $p < p_2$, where the dynamics, confined to the largest component, allows for a direct measurement of the impact of low connectivity on spreading time. Notably, for the reference parameter configuration with $N=10^{4}$ and $Z_0=1$, the discrete transient times reach a peak at $p=p_1= 10^{-4}$, which aligns with the classical critical threshold. Analysing this phenomenon is essential as it reflects the interplay of two opposing factors affecting these times: firstly, the decrease in connectivity, which, as demonstrated earlier, decelerates the dynamics; and secondly, the reduced size of the giant component, which naturally accelerates the time required for complete infection. It is noteworthy to observe that the number of sizes or classes of components, as depicted in Fig. \ref{NGdep}b, also reaches a maximum at this critical point. Note that, in contrast to transient times, both the characteristic time and the reciprocal of the eigenvalue are decreasing functions of $p$.


\section{Discussion}

In this paper, we have studied the population dynamics of a Susceptible-Infected (SI) model on a random Erd\"os-R\'enyi network, comparing it with the mean field approach in terms of the parameter $p$. This mean field approach, expressed through Ordinary Differential Equations (ODEs), assumes a law of mass action for local contagion at network nodes, introducing the parameter $\beta$, which represents a rate of effectiveness of the contact between a susceptible individual with its infected neighbours \cite{Wilson, Voit}. Notably, $\beta$ is absent in the discrete network model, where a probabilistic contagion rule is considered.

To demonstrate the relationship between the macroscopic contagion rate $\beta$, connectivity, and contagion spread, as already stated by some authors (see, for instance, \cite{Kiss, Aleta}), we calibrated $\beta$ in terms of the linking probability $p$ of the Erd\"os-R\'enyi random network. We identified a critical $p$-value, $p_d$, such that for $p > p_d$ the mean field approach is valid. This implies that the time evolution of the average curve aligns between the simulations and the analytical solution of the Initial Value Problem associated to the ODE \ref{eqZ} with $\beta = 1$, independently of the initial infected population. Interestingly, for $p >p_d$, changes in average connectivity have no effect on the dynamics, which could be relevant to implement measures against an epidemic that follows this SI-model.

Just the contrary occurs for lower $p$-values as shown in this paper. Indeed,  for $p <p_d$ the mean field approach loses validity due to (i) the dependence on the initial infected population and (ii) the value of $\beta$ less than 1. To analyse the dependence of $\beta$ on the average connectivity of the network $k$ or, equivalently, with $p$, we constrained the dynamics to the largest component for each $p$,  eliminating the effects of network fragmentation at low $p$-values. The obtained result indicates that $\beta$ is a monotonically increasing function of $p$ (see Fig.{panel}a), ranging from 0 at $p=0$ to 1 at $p=1$. Additionally, the averaged standard error for $\beta$ peaks at the classical threshold $p_1$, as shown in Fig. \ref{panel}b.

To validate these findings, we computed the time the infected population takes to reach equilibrium using both simulations (transient times $t_{\rho}$) and the macroscopic ODEs (characteristic time $T_c$). Consistent with the $\beta$ estimation, for $p >p_d$, these times remain practically constant. Conversely, for $p < p_d$, the transient times exhibit a maximum at $p=p_1$ (see Fig.\ref{tt}), revealing a previously unreported coincidence.

\section{Further Considerations}

In this paper, we have assumed that the classical law of mass action applies, and it is the contagion rate $\beta$ which is calibrated as a function of the network connectivity. However, generalizations of this law exist that consider a different functional macroscopic response to the microscopic contagion process that occur at each node of the network. For instance, it has been postulated that more complex interactions between susceptible and infected individuals could be caught by using power-laws as: $\beta \, X^a \, Z^b$, where $a$ and $b$ are exponents to be experimentally determined and $\beta$ can be dependent on network structure as well \cite{Keeling}.  This kind of models should be typically applied to complex contagion process, e.g. other types of social transmission. Nonetheless, for the simple contagion SI model considered in this paper, the simplest assumption of $a=b=1$ seems more suitable, taking $\beta$ a function of $p$. 

Another avenue for exploration involves alternative generalizations that account for non-homogeneous and non-mixing populations. Models based on connected pairs in the network, represented through higher-order ODE systems, provide a framework for such situations \cite{Kiss, Keeling}.  Comparing the results obtained from these models, assuming the classical law of mass action with a $\beta$ dependent on $p$, could offer valuable insights.


%

\end{document}